\documentclass[a4paper,12pt]{amsart}

\usepackage{amssymb}
\usepackage{dsfont}
\usepackage{enumerate}
\usepackage[inline]{enumitem}
\usepackage{mathtools}
\usepackage{amsmath}
\usepackage{color}
\usepackage{bbm}         
\usepackage{amsfonts}
\usepackage{graphicx,hyperref,verbatim}
\usepackage{amsthm}
\usepackage{thmtools}
\usepackage{cleveref}
\usepackage{url}
\usepackage{tikz,multicol,manfnt}
\usepackage{qtree}
\usepackage[weather,clock,alpine]{ifsym}
\usepackage{xparse}
\usepackage{parnotes}
\usepackage{placeins}
\usepackage{todonotes}

\newcommand{\rond}{\mathcal}

\NewDocumentCommand{\set}{mg}{\left\{#1\IfNoValueF{#2}{\;\middle|\;#2}\right\}}

\newcommand{\ignore}[1]{}

\let\phi\varphi

\DeclareMathOperator{\rk}{rk}       
\DeclareMathOperator{\drk}{drk}

\theoremstyle{plain}
\newtheorem{thm}{Theorem}[section]
\newtheorem{cor}[thm]{Corollary}
\newtheorem{lem}[thm]{Lemma}
\newtheorem{prop}[thm]{Proposition}

\newtheorem{fct}[thm]{Fact}
\newtheorem{ans}[thm]{Answer}
\theoremstyle{definition}
\newtheorem{dfn}[thm]{Definition}
\theoremstyle{remark}

\newtheorem{qst}[thm]{Question}
\newtheorem{exam}[thm]{Example}

\crefname{exam}{example}{examples}

\todostyle{blaise}{author=blaise,color=pink}
\todostyle{lasse}{author=lasse,color=cyan}

\title{Know your Rank!}
\author{Blaise Boissonneau}
\author{Lasse Vogel}
\date{\today}

\begin{document}

\begin{abstract}
We study definable ranks of ordered fields, ordered abelian groups, and linear orders. For an arbitrary linear order $\Gamma$, we construct an ordered abelian group $G$ with archimedian spine $\Gamma$ and an ordered field $K$ with natural value group $G$ such that the definable ranks of $K$, $G$ and $\Gamma$ are all isomorphic. This answers a question of Krapp, Kuhlmann, and the second author.
\end{abstract}

\maketitle


\section{Introduction}
This papers focuses on several questions raised in \cite{KKV}, as such, it is relatively short. It builds on discussions which took place during the conference ``Structures algébriques et ordonnées'' organized in the occasion of the 40th anniversary of the DDG seminar. Despite being thought as a companion of the aforementionned paper, it is self-contained. We repeat here relevant definitions and all the results we need.

\begin{dfn} Let $(K,<)$ be an ordered field. The \emph{rank} of $(K,<)$, denoted $\rk(K)$, is the set of all convex non-trivial valuation rings of $K$, ordered by inclusion.

We define an equivalence relation saying that elements belong to the same archimedian class, specifically, $a\sim b$ iff $\exists n\in\mathds N$ such that $|a|<n|b|$ and $|b|<n|a|$ (where $|\cdot|$ is the absolute value). Let $[\cdot]$ denote the equivalence class modulo $\sim$. We will write $[a]>[b]$ if $|a|<|b|$ and $a\not\sim b$ (note the inversion) and $[a]+[b]$ will denote $[ab]$. We now call \emph{natural valuation} of $K$ the canonical projection $v_{nat}\colon K\rightarrow K/{\sim}$.

Finally, the \emph{definable rank} of $(K,<)$, denoted $\drk(K)$, is the set of all definable convex non-trivial valuation rings of $K$, ordered by inclusion.
\end{dfn}

\begin{dfn}
Let $(G,<)$ be an ordered abelian group. The \emph{rank} of $(G,<)$, denoted $\rk(G)$, is the set of all convex proper subgroups of $G$, ordered by inclusion.

We define, similarly as above, an equivalence relation saying that elements belong to the same archimedian class, and we order the archimedian classes as above. We now call \emph{natural valuation} the canonical projection $v_G\colon G\mapsto G/{\sim}$. The image $v_G(G)$ is a linearly ordered set. We will call $v_G(G\setminus\{0\})$ the \emph{archimedian spine} of $G$.

Finally, the \emph{definable rank} of $(G,<)$, denoted $\drk(G)$, is the set of all definable convex proper subgroups of $G$, ordered by inclusion.
\end{dfn}

\begin{dfn}
Let $(\Gamma,<)$ be a linear order. The \emph{rank} of $(\Gamma,<)$, denoted $\rk(\Gamma)$, is the set of all proper end segments of $\Gamma$, ordered by inclusion.

The \emph{definable rank} of $(\Gamma,<)$, denoted $\drk(\Gamma)$, is the set of all definable proper end segments of $\Gamma$, ordered by inclusion.
\end{dfn}

The aforementionned paper states the following:

\begin{fct} [{\cite[Fact~3.3]{KKV}}]
 Let $(K,<)$ be an ordered field, $G=v_{nat}(K^\times)$ its natural value group, $\Gamma=v_G(G\setminus\{0\})$ its archimedian spine. Then there are natural isomorphisms between the ranks:
 
 $$\rk(K)\xrightarrow{\Phi_K}\rk(G)\xrightarrow{\Phi_G}\rk(\Gamma)$$

 defined by $\Phi_K(\rond O)=v_{nat}(\rond O^\times)$ and $\Phi_G(H)=v_G(H\setminus\{0\})$.
\end{fct}

So the three notions of rank are in exact correspondance. A natural question is to ask how much of this remains true when considering definable ranks. In \cite{KKV} the authors give examples where all three definable ranks are pairwise non-isomorphic. We will in turn discuss some such examples, see \cref{drkG≠drkΓ,drkK=drkG≠drkΓ} but also and more importantely, examples where definable ranks are indeed isomorphic, see \cref{drkK=drkG=drkΓ}.

For our discussion today, we need results guaranteeing isomorphicity of definable ranks. Namely, we repeat the following:

\begin{prop}[{\cite[Lemma~5.5]{KKV}}] Let $(K,<)$ be an ordered field and let $G=v_{nat}(K^\times)$. If $v_{nat}$ is henselian, then $\Phi_K(\drk(K))$ is an end segment of $\drk(G)$.
\end{prop}

It follows that $\Phi_K(\drk(K))=\drk(G)$ if and only if the convex subgroup $\{0\}$ belongs to $\Phi_K(\drk(K))$ if and only if the valuation ring $\mathcal O_{nat}$ belongs to $\drk(K)$, that is:

\begin{cor}\label{vnatdef} Let $(K,<)$ be an ordered field with henselian natural valuation and $G=v_{nat}(K^\times)$. Then $\Phi_K(\drk(K))=\drk(G)$ if and only if $v_{nat}$ is definable.\footnotemark
\end{cor}

\footnotetext{This result was stated in a previous version of \cite{KKV}. Since then, the authors improved their result, mostly by classifying when is a henselian $v_{nat}$ definable in the first place. We only rely on the result which we sate here.}

We are now ready to attack the matter of the day. We will focus our attention on three questions asked in \cite{KKV}:



\begin{qst}[{\cite[Question 7.3]{KKV}}]\label{Q3} Given an ordered set $(C, <)$, when is it possible to construct an ordered field $K$ such that $\drk(\Gamma)\simeq(C, <)$, $\Phi_G^{-1}(\drk(\Gamma)) = \drk(G)$ and $\Phi_K^{-1}(\drk(G)) = \drk(K)$ where $G = v_{nat}(K^\times)$ and $\Gamma = v_G(G\setminus\set{0})$?
\end{qst}

\begin{qst}[{\cite[Question 7.4]{KKV}}]\label{Q4} Given an ordered set $(\Gamma, <)$, when is it possible to construct an ordered field $K$ such that $v_G(G\setminus\set{0}) = \Gamma$, $\Phi_K (\drk(K) ) = \drk(G)$  and $\Phi_G(\drk(G)) =\drk(\Gamma)$, where $G = v_{nat}(K^\times)$?
\end{qst}



\begin{qst}\label{Q6} Let $(\Gamma, <)$ be an ordered set. What are necessary conditions for an end segment $\Delta\subseteq\Gamma$ to be definable?\footnote{This question too was present in an older version of \cite{KKV} but has disapeared in the latest version.}
\end{qst}

We answer \cref{Q3} partially and \cref{Q4} completely by using the following construction:

\begin{ans}
Let $(\Gamma,<)$ be a linear order. Let $G=\sum_{\Gamma}\mathds Z$ and $K=\mathds R((G))$. Then $v_{nat}=(K^\times)=G$, $v_G(G\setminus\{0\})=\Gamma$, $\Phi_K(\drk(K))=\drk(G)$ and $\Phi_G(\drk(G)) =\drk(\Gamma)$.
\end{ans}

This construction proves that the answer to \cref{Q4} is ``always'', while  answering \cref{Q3} is dependent on classifying which linear orders $(C,<)$ can be the definable ranks of a linear order. Note that such a classification is, more or less, an answer to \cref{Q6}.

We will now present the construction. Afterwards, we will discuss some examples showing how complicated an answer to \cref{Q6} really is.

\section{Spines}

In order to study the definable rank of ordered abelian groups, we recall the formalism of spines developped by Gurevich and Schmitt. This is collected from Schmitt's habilitation, see \cite{Schmitt}.

We encourage readers to skip directly to \cref{qespines}, as it is the only relevent result for the rest of our discussion today.

\begin{dfn}
 Let $G$ be an ordered abelian group and let $g\in G$. We write ``$C\lhd G$'' to mean that $C$ is a convex subgroup of $G$ ($C=\set0$ and $C=G$ are allowed). We define the following convex subgroups:
 \begin{itemize}
 \item $A(g)=\bigcup\set{C\lhd G}{g\notin C}$ is the largest convex subgroup not containing $g$,
 \item $B(g)=\bigcap\set{C\lhd G}{g\in C}$ is the smallest convex subgroup containing $g$.
\end{itemize}

We immediately remark that $g$ and $h$ are in the same archimedian class iff $A(g)=A(h)$, so $v_G(G\setminus\set0)$ can equivalentely be seen as the familly $\set{A(g)}{g\in G\setminus\set0}$.

Let $n\geqslant2$. We say that an ordered abelian group $H$ is $n$-regular if any interval containing at least $n$ points contains an $n$-divisible point. We finally define:
\begin{itemize}
 \item $A_n(g)=\bigcap\set{C\lhd A(g)}{B(g)/C\text{ is }n\text{-regular}}$, making the family $\set{A_n(g)}{g\in G\setminus\{0\}}$ the ``$n$-regular rank'' of $G$, and
 \item $F_n(g)=\bigcup\set{C\lhd G}{C\cap(g+nG)=\emptyset}$, the $n$-fundament of $g$, which is the largest convex subgroup modulo which $g$ is not $n$-divisible.
 \end{itemize}

The definition of $A_n(g)$ runs into an edge case when $g=0$; in this case, and this case only, we will take $A_n(g)=\emptyset$, which is not a convex subgroup.

Similarly, if $g$ is $n$-divisble, we take $F_n(g)=\emptyset$.

For notational ease, we also define $E_n(g)=\set{h\in G}{F_n(h)\subseteq F_n(g)}$, $E_n^*(g)=\set{h\in G}{F_n(h)\subsetneq F_n(g)}$, and $F_n^*(g)=E_n(g)/E_n^*(g)$.
\end{dfn}

\begin{fct}\label{defspines}
For every $n\geqslant 2$, the $A_n$ and the $F_n$ are uniformly definable convex subgroups of $G$. Consequently, the set $$Sp_n(G)=\set{A_n(g)}{g\in G\setminus\{0\}}\cup\set{F_n(g)}{g\in G\setminus nG}$$ is interpretable in $G$.
\end{fct}

We equip $Sp_n(G)$ with the ordering $\subseteq$ and with various colours; 
namely:

\begin{dfn}
The $n$-spine of $G$ is the set $Sp_n(G)$ equipped with:
\begin{itemize}
\item a binary relation symbol $<$, to be interpreted as $\subsetneq$,
\item unary predicates $A$ and $F$, to be interpreted as $A(x)$ iff $x=A_n(g)$ for some $g\in G$ and similarly for $F$,
\item a unary predicate $D$, to be interpreted as $D(x)$ iff $G/x$ is discrete (that is, there exists a minimum positive element in $G/x$),
\item unary predicates $\alpha_{p,k,m}$ for $p$ prime dividing $n$, $k$ integer such that $0\leqslant k<v_p(n)$, and $m$ an arbitrary integer, to be interpreted as $\alpha_{p,k,m}(x)$ iff $\alpha_{p,k}(F_n^*(g))\geqslant m$, where $\alpha$ is the Szmielev-invariant, that is,
$$\alpha_{p,k}(H)=\dim_{\mathds F_p}((p^kH)[p]/(p^{k+1}H)[p])$$
where $K[p]$ is the $p$-torsion elements of a group $K$.
\end{itemize}
\end{dfn}

For us, these colours will mostly not matter. Thus, we recommend not to spend too much time dissecting the above definition, or the below definition.

\begin{dfn}
The language $\rond L_{OG}^3$ is a multisorted language composed of:
\begin{itemize}
\item one sort $G$ for the ordered abelian group, equipped with:
\begin{itemize}
\item the language of ordered abelian groups $\set{0,+,-,<}$,
\item unary predicates $M_{n,k}$ for $k>0$ and $n\geqslant 0$, to be interpreted as $M_{n,k}(g)$ iff $G/A_n(g)$ is discrete and $g+A_n(g)$ is equal to $k$ times the minimal positive element,
\item unary predicates $D_{p,r,i}$ for each prime $p$ and arbitrary integers $r,i$, to be interpreted as $D_{p,r,i}(g)$ iff $g$ is $p^r$-divisible or $g+E^*_{p^r}(g)$ is $p^i$-divisible in $F^*_{p^r}(g)$.
\item unary predicates $E_{p,r,k}$ for each prime $p$ and arbitrary integers $r,k$, to be interpreted as $E_{p,r,k}(g)$ iff there is $h\in G$ with $F_{p^r}(g)=A_{p^r}(h)$, $G/A_{p^r}(h)$ is discrete, $h+A_{p^r}(h)$ is its minimum positive element, and $F_{p^r}(g-kh)\subsetneq F_{p^r}(g)$;
\end{itemize}
\item and for each $n$ a sort $S_n$ for the $n$-spine, equipped as above with the relation $<$ and predicates $A$, $F$, $D$ (which we hate), and $\alpha_{p,k,m}$ for $p^{k+1}$ dividing $n$ and $m$ arbitrary.
\end{itemize}
Furthermore, we add to this language function symbols $A_n$ and $F_n$ for each $n$, coming from the main sort $G$ and going to the sort $S_n$, with the obvious interpretation. 
\end{dfn}

Note that this language is a definable/interpretable expansion of the language of ordered abelian groups.

The sole purpose of defining this language is the following:

\begin{fct}\label{qespines}
In this language, the theory of ordered abelian groups eliminates quantifiers from the main sort $G$. Consequently, as the sorts $S_n$ are closed, each spine $Sp_n(G)$ is ``pure with control over parameters'' (see \cite[Section~1.1.1]{Tou}), and in particular, any subset of $Sp_n(G)$ which is definable in the full language of ordered abelian groups is actually definable in the pure coloured linear order $Sp_n(G)$.
\end{fct}

\section{The construction}

We are now ready to construct an ordered abelian group and an ordered field with exactly the definable rank we want. We will first construct the group, and will construct it as a lexicographic sum, that is:

\begin{dfn}
Let $(\Gamma,<)$ be a linear order. Let $(G_i)_{i\in\Gamma}$ be a sequence of ordered abelian groups. Then we define $G=\sum_{i\in\Gamma}G_i$ as follows: as a set, it is the subset of the cartesian product of the $G_i$ containing only sequences with finite support, that is, $g\in G$ iff $\set{i\in\Gamma}{g_i\neq0}$ is finite. It is equipped with the coordinatewise group law, and is ordered lexicographically: $g>h$ iff $\exists i\in\Gamma$ such that $g_i>h_i$ and $g_j=h_j$ for all $j<i$. 
\end{dfn}

\begin{prop}
Let $(\Gamma,<)$ be a linear order and let $G=\sum_{\Gamma}\mathds Z$. Then $v_G(G\setminus\set0)\simeq\Gamma$, and (modulo this identification) $\Phi_G(\drk(G))=\drk(\Gamma)$.
\end{prop}

\begin{proof}
First, note that convex subgroups of $G$ are exactly of the form $\sum_{\Delta}\mathds Z$ where $\Delta$ is an end segment of $\Gamma$ -- as a convention, a sum indexed by the emptyset will equal $\set0$.

Let $g\in G\setminus\set0$. Let $i\in\Gamma$ be the largest coordinate of $g$ which is not $0$. Then for all $n\geqslant 2$, we have $A(g)=A_n(g)=\sum_{\Gamma_{>i}}\mathds Z$. We can thus identify $\Gamma$ and $v_G(G\setminus\set0)$ by the isomorphism $i\mapsto \sum_{\Gamma_{>i}}\mathds Z$.

Let $g\in G\setminus nG$. Let $i\in\Gamma$ be the largest coordinate of $g$ which is not in $n\mathds Z$. Then $F_n(g)=\sum_{\Gamma_{>i}}\mathds Z$.

This means that $(Sp_n(G),\subseteq)$ is isomorphic to $(\Gamma,<^{inv})$ (note the inversion) via $\Gamma\rightarrow Sp_n(G)$, $i\mapsto \sum_{\Gamma_{>i}}\mathds Z$. Computing the colours in $Sp_n$, we quickly see that all the points have the same colours, so we might simply consider $Sp_n(G)$ as a linear order.

Consider a non-trivial end segment $\Delta$ definable in $(\Gamma,<)$. The corresponding inital segment $\nabla$ of $(Sp_n(G),\subseteq)$ is also definable. Now, ``$F_n(g)\in\nabla$'' is a first-order formula on $g$, so
$$\sum_{\Delta}\mathds Z=\bigcup\set{F_n(g)}{g\in G,F_n(g)\in\nabla}$$
is definable.

Consider now a convex subgroup $H=\sum_{\Delta}\mathds Z$ definable in $G$. Then $\set{F_n(g)}{F_n(g)\subseteq H}$ is definable too. This means that the initial segment $\nabla$ of $Sp_n(G)$ corresponding to $\Delta$ is definable. By purity (see \cref{qespines}), this is now a definable set of $Sp_n(G)$ which is isomorphic to $\Gamma^{inv}$, so the end segment $\Delta$ is definable in $(\Gamma,<)$.
\end{proof}

\begin{exam}\label{drkG≠drkΓ}
 let $\Gamma$ be a linear order and $\Delta\subseteq\Gamma$ an end segment which is \emph{not} definable. Then $G=\sum_{\Gamma\setminus\Delta}\mathds Z_{(2)}\oplus\sum_{\Delta}\mathds Z_{(3)}$ has a strictly bigger definable rank than $\drk(\Gamma)$, in the sense that $\Phi^{-1}_G(\drk(\Gamma))$ is a strict subset of $\drk(G)$.

 Let $\Gamma$ be a linear order and $G=\sum_{\Gamma}\mathds Q$. Then $\drk(G)$ is strictly smaller than $\drk(\Gamma)$, in the sense that $\Phi_G(\drk(G))$ is a strict subset of $\drk(\Gamma)$.

 Playing with these two examples, we can easily construct an ordered abelian group such that $\drk(G)$ and $\drk(\Gamma)$ do not embedd in one another via $\Phi_G$.
\end{exam}

We now construct an ordered field with natural value group as above and with exactly the same definable rank. We will construct it as a Hahn Field:

\begin{dfn}
 Let $k$ be an ordered field and $G$ an ordered abelian group. We construct the ordered field $k((G))$ as the field of formal power series $\set{\sum_{g\in G}a_g t^g}{a_g\in K}$ with well-ordered support, that is, $\set{g\in G}{a_g\neq 0}$ is a well-ordered subset of $G$. We order $k((G))$ by saying that $t$ is infinitesimal positive, that is, $\sum_{g\in G}a_g t^g>0$ iff $a_{g_0}>0$, where $g_0=\min(\set{g\in G}{a_g\neq 0})$.

 In such a field, the valuation defined by $v(\sum_{g\in G}a_g t^g)=g_0$, where $g_0=\min(\set{g\in G}{a_g\neq 0})$ as above, is a henselian valuation with value group $G$ and residue field $k$. If $k$ is archimedian as an ordered field, then this valuation is the natural valuation $v_{nat}$ on $k((G))$.
\end{dfn}

\begin{prop}\label{drkK=drkG=drkΓ}
 Let $\Gamma$ be a linear order and $G=\sum_\Gamma\mathds Z$ as above. Then $K=\mathds R((G))$ has $v_{nat}(K^\times)=G$ and $\Phi_K(\drk(K))=\drk(G)$ (and thus also $\drk(\Gamma)$).
\end{prop}

\begin{proof}
It is clear that $v_{nat}(K^\times)=G$ since $\mathds R$ is archimedian. $v_{nat}$ is also henselian by what's above. By \cref{vnatdef}, we just have to prove that $v_{nat}$ is definable.

We use the method developped in \cite{itme-spines}, and particularly, the following fact:

\begin{fct}[{\cite[Theorem 2.21]{itme-spines}}]
 Let $(K,v)$ be a henselian valued field of equicharacteristic 0. Then $v$ is definable if and only if the residue field $Kv$ and value group $vK$ are not ``co-augmentable'', that is, there is no non-trivial ordered abelian group $\Delta$ such that both $Kv\preccurlyeq Kv((\Delta))$ and $vk\preccurlyeq vK\oplus\Delta$ hold.
\end{fct}

In our case, such a $\Delta$ would need to be such that $\mathds R\preccurlyeq\mathds R((\Delta))$. By classical results on real-closed fields, this forces $\Delta$ to be divisible.

Now $G$ (provided $\Gamma$ is non-empty) satisfies $\forall g(g\neq0\rightarrow\exists h (|h|\leqslant|g|\wedge h\notin nG))$, for $n\geqslant 2$. This means any $\Delta$ with $G\preccurlyeq G\oplus\Delta$ would also have this property, in particular, wouldn't be divisible.

Thus $G$ and $\mathds R$ aren't co-augmentable, and $v_{nat}$ is definable.
\end{proof}

This gives an example where the definable ranks are preserved between the field, the value group, and the archimedian spine. We previously gave example where the definable ranks differ betweens the value group and the archimedian spine. We now give examples proving that every situation can arise:

\begin{exam}\label{drkK=drkG≠drkΓ}
 Let $\Gamma$ be a linear order and $\Delta\subseteq\Gamma$ an end segment \emph{not} definable. Let $G=\sum_{\Gamma\setminus\Delta}\mathds Z_{(2)}\oplus\sum_{\Delta}\mathds Z_{(3)}$ and $K=\mathds R((G))$. Then $\drk(K)=\drk(G)\neq\drk(\Gamma)$.


 Let $K=\mathds R((\mathds Q))$, $G=v_{nat}(K^\times)=\mathds Q$ and $\Gamma=v_G(G\setminus\{0\})$. Then $\drk(K)\neq\drk(G)=\drk(\Gamma)$.

 Let $K=((\mathds Q\oplus\mathds Q))$, $G=v_{nat}(K^\times)=\mathds Q\oplus\mathds Q$ and $\Gamma=v_G(G\setminus\{0\})$. Then $\drk(K)\neq\drk(G)\neq\drk(\Gamma)$.
\end{exam}

\section{An answer, a partial answer, and a non-answer}

\begin{cor}[answer to \cref{Q4}]
 Given an ordered set $(\Gamma, <)$, it is always possible to construct an ordered field $K$ such that $v_G(G\setminus\set{0}) = \Gamma$, $\Phi_K (\drk(K) ) = \drk(G)$  and $\Phi_G(\drk(G)) =\drk(\Gamma)$, where $G = v_{nat}(K^\times)$. Namely, $K=\mathds R((\sum_{\Gamma}\mathds Z))$.
\end{cor}

\begin{cor}[partial answer to \cref{Q3}]
 Given an ordered set $(C, <)$. Assume that there is a linear order $(\Gamma,<)$ such that $\drk(\Gamma)\simeq(C, <)$. Then it is possible to construct an ordered field $K$ such that $\drk(\Gamma)\simeq(C, <)$, $\Phi_G^{-1}(\drk(\Gamma)) = \drk(G)$ and $\Phi_K^{-1}(\drk(G)) = \drk(K)$ where $G = v_{nat}(K^\times)$ and $\Gamma = v_G(G\setminus\set{0})$. Namely, $K=\mathds R((\sum_{\Gamma}\mathds Z))$.
\end{cor}

A total answer to \cref{Q3} would involve completely characterising which linear orders $(C,<)$ could be the definable rank of a linear order $(\Gamma,<)$, in particular, would give an answer to \cref{Q6}. In some sense, the following is an answer:

\begin{fct}[{\cite[Proposition~2.27]{itme-spines}}]
 An end segment $I$ of a linear order $C$ is definable if and only if the pair $(C\setminus I, I)$ is not ``co-augmentable'', that is, there is no non-empty linear order $X$ such that $C\setminus I\preccurlyeq C\setminus I+X$ and $X\preccurlyeq X+I$.
\end{fct}

However, despite being a complete characterization of definable end segments of linear orders, it is a non-answer: given a specific end segment of a specific linear order, proving non-co-\hspace{0pt}augmentability can be as hard as proving definability directly, as we will see in \cref{cursed}.

Regarding the more specific question of determining whether a given linear order $C$ is the definable rank of a linear order $\Gamma$, we give below a condition which any such $C$ has to fulfill, but it is far from sufficient.

\begin{lem}
Let $(\Gamma,<)$ be a linear order and let $C=\drk(\Gamma)$. Then $C$ does not have a dense part; more specifically, for any $a<b$ in $C$, there is $d\in [a,b)$ such that $d$ has a successor, and similarly there is $d'\in (a,b]$ such that $d'$ has a predecessor.
\end{lem}

\begin{proof}
 Let $a<b$ in $C$. By definition, $a$ and $b$ are non-empty end segments of $\Gamma$, and $\exists \gamma\in a\setminus b$. For any such $\gamma$, let $d=\Gamma_{\geqslant \gamma}=\set{x\in\Gamma}{x\geqslant \gamma}$ and $d'=\Gamma_{>\gamma}$. Both $d$ and $d'$ are definable end segments of $\Gamma$, $a\leqslant d<d'\leqslant b$, and $d'$ is the successor of $d$.
\end{proof}

To see that this condition is not sufficient, consider $C=\sum_{\mathds Q} 3$; it satisfies said condition, but is not the definable rank of a linear order.

We finally discuss some examples of definable and non-definable end segments.

\begin{exam}\label{cursed}
Given a sequence $a\in\mathds N^\omega$, we consider the linear order $\Gamma_a=\sum_{k
\in\omega}(\mathds Q+a_k+2)$. It is not hard to see that $(\Gamma_a,<)$ interprets $a$.

Consider $\Gamma_{0,0,\dots}+\Gamma_{0,0,\dots}^{inv}$. The (end segment corresponding to the) cut in the middle is not definable, since the pair $(\Gamma_{0,0,\dots},\Gamma_{0,0,\dots}^{inv})$ is co-augmentable by $\sum_{k\in\mathds Z}\mathds Q+2$.

Consider $\Gamma_{0,0,\dots}+\Gamma_{1,1,\dots}^{inv}$. The cut in the middle is now definable, since $\Gamma_{1,1,\dots}^{inv}$ is the right-closure of the set of discrete sequences of 3 points.

Now consider $\Gamma_\pi+\Gamma_e^{inv}$, where $\pi$ and $e$ are the sequences of digits of $\pi$ and $e$ in base ten. Is the cut in the middle definable? The answer would be yes if, for example, $\pi$ is normal but $e$ isn't, or the opposite. The answer is probably no, as there is no pattern in the digits of either of those numbers that would make them distinguishable. But this is an open question in number theory.

Alternatively, one can also consider $a$ the sequence of digits of $\pi$ in base two and $b$ the digits of $\pi$ in base three. Now in $\Gamma_a+\Gamma_b^{inv}$, the cut in the middle would be definable if, for example, $\pi$ had infinitely many 2 in base three, or if ``$11$'' appeared only finitely many times in $\pi$ is base 2 but infinitely many times in $\pi$ in base 3, etc. Again, this is an open question in number theory.
\end{exam}

With this encoding, it becomes clear that answering \cref{Q6} and thus also \cref{Q4} in full generality will be very hard.

\bibliographystyle{abbrv}
\bibliography{ref}

\end{document}